\documentclass[11pt]{article}

\usepackage{tikz}
\usetikzlibrary{calc}
\usetikzlibrary{matrix}
\usepackage{amsmath,amssymb,amsfonts} % Typical maths resource packages
\usepackage{verbatim}
\usepackage{hyperref}                 % For creating hyperlinks in cross references

\usepackage{graphicx}

\def \p{\partial}
\newcommand{\V}[1]{\mbox{\boldmath $ #1 $}}
\newcommand{\bey}{\begin{eqnarray}}
\newcommand{\eey}{\end{eqnarray}}
\newcommand{\nn}{\nonumber}
\newcommand{\beq}{\begin{equation}}
\newcommand{\eeq}{\end{equation}}
\newtheorem{thm}{\hspace{6mm}Theorem}[section]

\newtheorem{lem}{\hspace{6mm}Lemma}[section]

\newcommand{\proofend}{\mbox{ }\hfill \raisebox{.4ex}{\framebox[1ex]{}}}
\vspace{5mm}

\begin{document}

\date{}
\title{Discrete maximum principle and a Delaunay-type mesh condition
for linear finite element approximations of two-dimensional anisotropic diffusion problems
}
\author{Weizhang Huang \thanks{
Department of Mathematics, the University of Kansas, Lawrence, KS 66045,
U.S.A. ({\tt huang@math.ku.edu}).
The work was supported in part by the National Science Foundation (USA)
under Grant DMS-0712935.}
}
\maketitle

\vspace{10pt}

\begin{abstract}
The finite element solution of two-dimensional anisotropic diffusion problems is considered.
A Delaunay-type mesh condition is developed for linear finite element approximations to satisfy
a discrete maximum principle. The condition is shown to be weaker than the existing
anisotropic non-obtuse angle condition. It reduces to the well known Delaunay condition
for the special case with the identity diffusion matrix. Numerical results are presented to
verify the theoretical findings.
\end{abstract}



\vspace{10pt}

% section 1
\section{Introduction}

We are concerned with the linear finite element (FEM) solution of the two-dimensional anisotropic diffusion
equation
\beq
\label{bvp-pde}
 - \nabla \cdot (\mathbb{D} \, \nabla u)  =  f, \qquad \mbox{ in } \quad \Omega
\eeq
subject to the Dirichlet boundary condition
\beq
\label{bvp-bc}
u =  g, \qquad \mbox{ on } \quad \partial \Omega
\eeq
where $\Omega \in \mathbb{R}^2 $ is a connected polygonal domain,
$f$ and $g$ are given functions, and $\mathbb{D}=\mathbb{D}(x,y)$ is the diffusion matrix
assumed to be symmetric and strictly positive definite on $\Omega$.
This boundary value problem (BVP) is a model of anisotropic diffusion problems
arising in various fields such as plasma physics
\cite{GL09, GLT07, GYK05, NW00, SH07, Sti92},
petroleum reservoir simulation \cite{ABBM98a,ABBM98b,CSW95,EAK01,MD06},
and image processing \cite{CS00, CSV03, KM09, MS89, PM90, Wei98}.
A distinct feature of the BVP is that its solution satisfies the maximum principle and
is monotone when $f(x,y) \le 0$ for all $(x,y) \in \Omega$.
A challenge in the numerical solution of the BVP is to design a scheme so that
the resulting numerical approximations satisfy a discrete maximum principle (DMP).

Development of DMP satisfaction schemes for solving diffusion problems has attracted considerable
interest in the past; e.g., see
\cite{BKK08,BE04,Cia70,CR73,KK05,KK06,KKK07,KL95,Let92,Sto82,Sto86,SF73,Var66,XZ99}
for isotropic diffusion problems where $\mathbb{D} = a(x,y) I$
with $a(x,y)$ being a scaler function and
\cite{CSW95,DDS04,GL09,GLT07,GYK05,KSS09,LePot05,LePot09,LePot09b,LH10,LSS07, LSSV07,LS08,MD06,SH07}
for anisotropic diffusion problems where $\mathbb{D}(x,y)$ can be heterogeneous and anisotropic.
For example, Ciarlet and Raviart \cite{CR73} (also see Brandts et al. \cite{BKK08}) show
that the linear finite element method for an isotropic diffusion problem
satisfies DMP when the mesh is simplicial and satisfies the non-obtuse angle
condition requiring the dihedral angles of mesh elements to be non-obtuse.
In two dimensions and for the special case $\mathbb{D} = I$,
the condition can be replaced by the Delaunay condition, a weaker condition that only requires
the sum of any pair of angles opposite a common edge to be less than or equal to $\pi$ \cite{Let92,SF73}.
Moreover, Xu and Zikatanov \cite{XZ99} show that the non-obtuse angle condition at edges where
the diffusion coefficient is discontinuous and the Delaunay condition at other places guarantee
DMP satisfaction. Recently, Li and Huang \cite{LH10} generalize the non-obtuse angle condition
to anisotropic diffusion problems and obtain the so-called anisotropic non-obtuse angle condition 
requiring the dihedral angles of mesh elements to be non-obtuse
when measured in a metric depending on $\mathbb{D}$.

The objective of this paper is to extend the Delaunay condition to anisotropic diffusion problems.
A Delaunay-type mesh condition is developed for the DMP satisfaction of linear finite element approximations
for those problems. It is shown that the new condition reduces to the Delaunay condition for the special
case $\mathbb{D} = I$ and is weaker than the anisotropic non-obtuse angle condition developed in \cite{LH10}.
We attain the new condition by investigating the stiffness matrix as a whole. This is different from
\cite{LH10} where only local stiffness matrices on individual elements are considered.
The main theoretical result is given in Theorem \ref{thm4.1}.

This paper is organized as follows. The linear finite element formulation for BVP (\ref{bvp-pde})
and (\ref{bvp-bc}) is given in Section 2. Section 3 is devoted to the description and geometric interpretation
of the anisotropic non-obtuse angle condition. The Delaunay-type mesh condition is developed 
in Section 4, followed by Section 5 with numerical results verifying the theoretical findings.
Finally, Section 6 contains conclusions and comments.

% section 2
\section{Linear finite element formulation for the model problem}

Consider the linear finite element solution of BVP (\ref{bvp-pde}) and (\ref{bvp-bc}).
Assume that a family of triangular meshes $\{ \mathcal{T}_h \}$ is given for $\Omega$.
Let
\[
U_g = \{ v \in H^1(\Omega) \; | \: v|_{\p \Omega} = g\}.
\]
Denote by $U_{g^h}^h \subset U_g$ the linear finite element space associated with mesh $\mathcal{T}_h$,
where $g^h$ is a linear approximation to $g$ on the boundary.
A linear finite element solution $\tilde{u}^h \in U_{g^h}^h$ to BVP (\ref{bvp-pde})
and (\ref{bvp-bc}) is defined by
\beq
\label{disc-1}
\sum_{K \in \mathcal{T}_h} \int_{K} (\nabla v^h)^{T} \, \mathbb{D} \,
\nabla \tilde{u}^h d x d y =
 \sum_{K \in \mathcal{T}_h} \int_{K} f \, v^h d x d y, \quad \forall \; v^h \in U_0^h
\eeq
where $U_0^h = U_{g^h}^h$ with $g^h = 0$. Generally speaking, the integrals in (\ref{disc-1})
cannot be carried out analytically and numerical quadrature
is often necessary. We assume that a quadrature rule has been chosen on the reference element $\hat{K}$,
\beq
\label{quadra-rule}
\int_{\hat{K}} v(\xi, \eta) d \xi d \eta \approx |\hat{K}| \sum_{k=1}^m \hat{w}_k v(\hat{b}_k),
\quad \sum_{k=1}^m \hat{w}_k =1,
\eeq
where $\hat{w}_k$'s are the weights and $\hat{b}_k$'s are the quadrature nodes.
Many quadrature rules can be used for this purpose; e.g., see \cite{EG04}.
An example is $\hat{w}_k=\frac{1}{3}$ ($k=1,2,3$)
and the barycentric coordinates of the nodes
($\frac{1}{6},\frac{1}{6},\frac{2}{3}$), ($\frac{1}{6},\frac{2}{3},\frac{1}{6}$),
and ($\frac{2}{3},\frac{1}{6},\frac{1}{6}$).

Let $F_K$ be the affine mapping from $\hat{K}$ to $K$ such that $K = F_K(\hat{K})$, and denote
$b_k^K = F_K(\hat{b}_k)$, $k=1,\cdots,m$. Upon applying (\ref{quadra-rule}) to the integrals in (\ref{disc-1})
and changing variables, the finite element approximation problem becomes seeking $u^h \in U_{g^h}^h$ such that
\beq
\label{disc-2}
\sum_{K \in \mathcal{T}_h} |K| \sum_{k=1}^m \hat{w}_k
\; (\nabla v^h|_K)^T \; \mathbb{D}(b_k^K) \; \nabla u^h|_K
=  \sum_{K \in \mathcal{T}_h} |K| \sum_{k=1}^m \hat{w}_k f(b_k^K) \; v^h (b_k^K),
\quad \forall v^h \in U_0^h
\eeq
where $\nabla v^h|_K$ and $\nabla u^h|_K$ denote the restriction of $\nabla v^h$
and $\nabla u^h$ on $K$, respectively. We have used the fact that $\nabla v^h|_K$ and
$\nabla u^h|_K$ are constant in deriving (\ref{disc-2}). Let
\beq
\label{def-DK}
\mathbb{D}_K = \sum_{k=1}^m \hat{w}_k \mathbb{D}(b_k^K) .
\eeq
Obviously, $\mathbb{D}_K$ is an average of $\mathbb{D}$ on $K$.
Eq. (\ref{disc-2}) can be written into
\beq
\label{disc-3}
\sum_{K \in \mathcal{T}_h} |K| \; (\nabla v^h|_K)^T \; \mathbb{D}_K \; \nabla u^h|_K
=  \sum_{K \in \mathcal{T}_h} |K| \sum_{k=1}^m \hat{w}_k f(b_k^K) \; v^h(b_k^K),
\quad \forall v^h \in U_0^h.
\eeq

We now express (\ref{disc-3}) in a matrix form.
Denote the numbers of the elements, vertices, and interior vertices of mesh $\mathcal{T}_h$
by $N$, $N_v$, and $N_{vi}$, respectively. Assume that the vertices are ordered in such a way that
the first $N_{vi}$ vertices are the interior vertices. Then $U_0^h$ and $u^h$ can be expressed as
\beq
U_0^h = \text{span} \{ \phi_1, \cdots, \phi_{N_{vi}} \} ,
\eeq
\beq
\label{soln-approx}
u^h = \sum_{j=1}^{N_{vi}} u_j \phi_j + \sum_{j=N_{vi}+1}^{N_{v}} u_j \phi_j ,
\eeq
where $\phi_j$ is the linear basis function associated with the $j$-th vertex, $\V{a}_j$.
The boundary condition (\ref{bvp-bc}) is approximated by
\beq
u_j = g(\V{a}_j), \quad j = N_{vi}+1, ..., N_v .
\label{fem-bc}
\eeq
Substituting (\ref{soln-approx}) into and taking $v^h = \phi_i$ ($i=1, ..., N_{vi}$) in (\ref{disc-3})
and combining the resulting equations with (\ref{fem-bc}), we obtain the linear algebraic system
\beq
\label{fem-linsys}
A \, \V{u} = \V{f},
\eeq
where
\beq
A = \left [\begin{array}{cc} A_{11} & A_{12} \\ 0 & I \end{array} \right ],
\label{fem-matrix}
\eeq
\[
\V{u} = (u_1,..., u_{N_{vi}}, u_{N_{vi}+1},..., u_{N_v})^T,
\]
\[
\V{f} = ( f_1, ..., f_{N_{vi}}, g_{N_{vi}+1}, ..., g_{N_v} )^T,
\]
and $I$ in (\ref{fem-matrix}) is the identity matrix of size $(N_{v} - N_{vi})$.
The entries of the stiffness matrix $A$ and the right-hand-side
vector $\V{f}$ are given by
\beq
\label{stiffness-matrix}
a_{ij}= \sum_{K \in \mathcal{T}_h} |K| \;(\nabla \phi_i|_K)^T \; \mathbb{D}_K \; \nabla \phi_j|_K,
\quad i=1, ..., N_{vi},\; j=1, ..., N_{v}
\eeq
\beq
f_i = \sum_{K \in \mathcal{T}_h} |K| \sum_{k=1}^m \hat{w}_k f(b_k^K)\; \phi_i (b_k^K),
\quad i=1, ..., N_{vi}.
\eeq

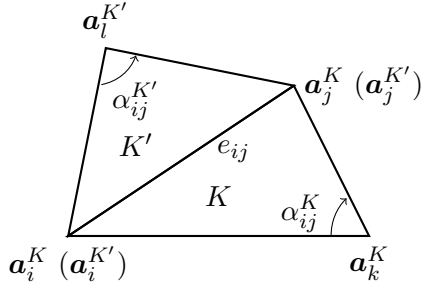
\begin{figure}[thb]
\centering
\begin{tikzpicture}[scale = 1]
\draw [thick] (-1,0) -- (3,0) -- (2, 2) -- cycle;
\draw [thick] (-1,0) -- (2, 2) -- (-0.5, 2.5) -- cycle;
\draw [below] (-1,0) node {$\V{a}_i^K$ ($\V{a}_i^{K'}$)};
\draw [below] (3,0) node {$\V{a}_k^K$};
\draw [right] (2,2) node {$\V{a}_j^K$ ($\V{a}_j^{K'}$)};
\draw [above] (-0.5,2.5) node {$\V{a}_l^{K'}$};
\draw [<-] (2.7,0.5) arc (135:180:0.7);
\draw [left] (2.5,0.3) node {$\alpha_{ij}^K$};
\draw [<-] (-0.1,2.4) arc (-20:-80:0.6);
\draw [below] (-0.1,2.2) node {$\alpha_{ij}^{K'}$};
\draw [below] (1,0.8) node {$K$};
\draw [below] (-0.1,1.5) node {$K'$};
\draw [below] (1.2,1.4) node {$e_{ij}$};
\end{tikzpicture}
\caption{Elements $K$ and $K'$ share the common edge $(e_{ij})$ connecting vertices
$\V{a}_i^K$ ($\V{a}_i^{K'}$) and $\V{a}_j^K$ ($\V{a}_j^{K'}$).
The angles opposite the edge are denoted by $\alpha_{ij}^K$ and $\alpha_{ij}^{K'}$, respectively. The Delaunay
condition is $\alpha_{ij}^K + \alpha_{ij}^{K'} \le \pi$. }
\label{f1}
\end{figure}

The expression (\ref{stiffness-matrix}) can be simplified. Let $\omega_i$ be
the patch of the elements sharing vertex $\V{a}_i$. Noticing that $\nabla \phi_i = 0$ for $(x,y) \notin \omega_i$,
we have, for $ i \neq j$, $i=1, ..., N_{vi}$, $j=1, ..., N_{v}$,
\bey
a_{ij} & = & \sum_{K \in \omega_i\cap \omega_j } |K| \;(\nabla \phi_i|_K)^T \; \mathbb{D}_K \; \nabla \phi_j|_K
\nn \\
& = & |K| \;(\nabla \phi_i|_K)^T \; \mathbb{D}_{K} \; \nabla \phi_j|_K
+ |K'| \;(\nabla \phi_i|_{K'})^T \; \mathbb{D}_{K'} \; \nabla \phi_j|_{K'} .
\label{stiffness-matrix-1}
\eey
 In (\ref{stiffness-matrix-1}), $K$ and $K'$ denote the two elements sharing the common edge
($e_{ij}$) connecting vertices $\V{a}_i \equiv \V{a}_i^K \equiv \V{a}_i^{K'}$ and $\V{a}_j \equiv
\V{a}_j^K\equiv \V{a}_j^{K'}$; see Fig. \ref{f1}.

% f1

% section 3
\section{The anisotropic non-obtuse angle condition}

In this section, we study mesh conditions under which the linear finite element scheme (\ref{disc-3})
satisfies DMP.

\begin{figure}[t]
\centering
\begin{tikzpicture}[scale = 1]
\draw [thick] (-1,0) -- (3,0) -- (2, 2) -- cycle;
\draw [below] (-1,0) node {$\V{a}_3^K$};
\draw [below] (3,0) node {$\V{a}_1^K$};
\draw [right] (2,2) node {$\V{a}_2^K$};
\draw [->] (-0.6,0) arc (0:30:0.42);
\draw [right] (-0.5,0.25) node {$\alpha_{12}^K$};
\draw [<-] (2.75,0.5) arc (135:180:0.7);
\draw [left] (2.5,0.3) node {$\alpha_{23}^K$};
\draw [->] (1,-0.4) -- (1,0.4);
\draw [below] (1,-0.4) node {$\V{q}_2^K$};
\draw [->] (3.0,1.25) -- (2,0.75);
\draw [right] (3.0,1.25) node {$\V{q}_3^K$};
\draw [->] (0.1,1.6) -- (0.82,0.52);
\draw [right] (0.1,1.6) node {$\V{q}_1^K$};
\draw [-] (3,0) -- (1.77,1.846);
\draw [left] (2.1,1.2) node {$h_1^K$};
%\draw [-] (-1,0) -- (2.2,1.6);
%\draw [left] (1.7,0.9) node {$h_3$};
\end{tikzpicture}
\caption{A sketch of the $\V{q}$ vectors and other geometric quantities for an arbitrary element $K$.}
\label{f2}
\end{figure}
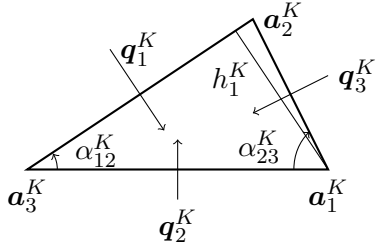

To start with, we introduce some notation.
Denote the vertices of  an element $K$ by $\V{a}_1^K, \V{a}_2^K, \V{a}_3^K$.
The edge matrix of $K$ is defined as
\beq
\nn
E_K = [\V{a}_2^K-\V{a}_1^K, \,\V{a}_3^K-\V{a}_1^K ].
\eeq
Since $K$ is simplicial, $E_K$ is nonsingular \cite{Som29}.
A set of $\V{q}$-vectors (cf. Fig. \ref{f2}) can then be defined as
\beq
\label{def-q}
[\V{q}_2^K, \, \V{q}_3^K] = E_K^{-T},
\quad \V{q}_1^K = - \V{q}_2^K - \V{q}_3^K.
\eeq
By definition, $\V{q}_i^K$ is the inward normal to the edge opposite to vertex $\V{a}_i^K$
(i.e., the edge not having $\V{a}_i^K$ as a vertex). This orthogonality implies that
the (dihedral) angle, $\alpha_{ij}^K$, opposite to edge $e_{ij}$
can be calculated in terms of $\V{q}_i^K$ and $\V{q}_j^K$ as
\beq
\label{dihedral}
\alpha_{ij}^K = \pi - \arccos\left ( \frac{\V{q}_i^K \cdot \V{q}_j^K }{ \|\V{q}_i^K\|\cdot \|\V{q}_j^K\| }\right ),
\quad i \ne j.
\eeq
Moreover, it is known 
\cite{BKK07,KL95} that
\beq
\label{q-phi}
\nabla \phi_{i}|_K = \V{q}_{i}^K.
\eeq
From this relation, it is not difficult to show
\beq
\| \V{q}_{i}^K \| = \frac{1}{h_i^K},
\label{q-h}
\eeq
where $h_i^K$ is the height of $K$ in the direction of $\V{q}_i^K$ or the shortest distance from
$\V{a}_i^K$ to the edge opposite to $\V{a}_i^K$; see Fig. \ref{f2}. 

% f2

Now, we are ready to describe the anisotropic non-obtuse angle condition.

\begin{lem}
\label{thm-g-nonobtuse}
If the mesh satisfies the anisotropic non-obtuse angle condition
\beq
\label{g-nonobtuse}
(\V{q}_i^K)^T \, \mathbb{D}_K \, \V{q}_j^K \le 0,
\quad \forall \; i \ne j,\; i,j = 1,2, 3,\; \forall \; K \in \mathcal{T}_h
\eeq
then the linear finite element scheme (\ref{disc-3}) for solving BVP (\ref{bvp-pde}) and (\ref{bvp-bc})
satisfies DMP.
\end{lem}

This lemma was proven in \cite{LH10} in any spatial dimension by showing
that the stiffness matrix $A$ in (\ref{fem-linsys}) is an $M$-matrix and has non-negative
row sums. A key step of the proof is to show $a_{ij} \le 0$ for all $i\neq j$,  which
can be seen to hold from  (\ref{stiffness-matrix-1}),  (\ref{q-phi}), and (\ref{g-nonobtuse}).

For the isotropic diffusion case,
the condition (\ref{g-nonobtuse}) reduces to 
\beq
\label{nonobtuse}
(\V{q}_i^K)^T \, \V{q}_j^K \le 0,
\quad \forall \; i \ne j,\; i,j = 1,2, 3,\; \forall \; K \in \mathcal{T}_h .
\eeq
Thus, (\ref{g-nonobtuse}) is a generalization of (\ref{nonobtuse})
for a general diffusion matrix.
Notice that (\ref{nonobtuse}) implies that the second angle on the right-hand side of
(\ref{dihedral}) is between $\pi/2$ and $\pi$.  Consequently, 
(\ref{nonobtuse}) is exactly the non-obtuse angle condition \cite{CR73},
implying $\alpha_{ij}^K \le \pi/2$.

The condition (\ref{g-nonobtuse}) can be more directly interpreted 
as requiring the angles of elements to be non-obtuse
when measured in a metric depending on $\mathbb{D}$.
To see this, we first notice that, according to (\ref{g-nonobtuse}),
the angle between $\V{q}_i^K$ and $\V{q}_j^K$
should be measured in the metric $\mathbb{D}_K$. Indeed, it has the expression
\[
\arccos \left (\frac{(\V{q}_i^K)^T \mathbb{D}_K \V{q}_j^K}
{\| \V{q}_i^K\|_{\mathbb{D}_K} \; \| \V{q}_j^K\|_{\mathbb{D}_K}} \right ) ,
\]
where the $\mathbb{D}_K$-norm is defined by
\beq
\| \V{v} \|_{\mathbb{D}_K} = \sqrt{\V{v}^T \mathbb{D}_K \V{v}},\quad \forall\; \V{v} \in \mathbb{R}^2 .
\label{M-norm}
\eeq
Since
\[
\arccos \left (\frac{(\V{q}_i^K)^T \mathbb{D}_K \V{q}_j^K}
{\| \V{q}_i^K\|_{\mathbb{D}_K} \; \| \V{q}_j^K\|_{\mathbb{D}_K}} \right )
= \arccos \left (\frac{(\mathbb{D}_K^{\frac{1}{2}} \V{q}_i^K)^T (\mathbb{D}_K^{\frac{1}{2}} \V{q}_j^K)}
{\| \mathbb{D}_K^{\frac{1}{2}} \V{q}_i^K\| \cdot \| \mathbb{D}_K^{\frac{1}{2}} \V{q}_j^K \|}\right ),
\]
the angle can also be regarded as the one between vectors $\mathbb{D}_K^{\frac{1}{2}} \V{q}_i^K$ and
$\mathbb{D}_K^{\frac{1}{2}} \V{q}_j^K$ in the Euclidean norm. Denote the third vertex of $K$ by
$\V{a}_k^K$. By definition, $\V{q}_i^K$ and $\V{q}_j^K$ are orthogonal to edges ($\V{a}_j^K - \V{a}_k^K$) and
($\V{a}_i^K - \V{a}_k^K$), respectively; i.e.,
\[
(\V{q}_i^K)^T (\V{a}_j^K - \V{a}_k^K) = 0,\quad
(\V{q}_j^K)^T (\V{a}_i^K - \V{a}_k^K) = 0 .
\]
It follows that
\[
(\mathbb{D}_K^{\frac{1}{2}}\V{q}_i^K)^T\left ( \mathbb{D}_K^{-\frac{1}{2}}(\V{a}_j^K - \V{a}_k^K) \right ) = 0,\quad
(\mathbb{D}_K^{\frac{1}{2}}\V{q}_j^K)^T \left ( \mathbb{D}_K^{-\frac{1}{2}}(\V{a}_i^K - \V{a}_k^K) \right )= 0 ,
\]
indicating that $\mathbb{D}_K^{\frac{1}{2}} \V{q}_i^K$ and $\mathbb{D}_K^{\frac{1}{2}}\V{q}_j^K$
are orthogonal to $\mathbb{D}_K^{-\frac{1}{2}}(\V{a}_j^K - \V{a}_k^K)$ and
$\mathbb{D}_K^{-\frac{1}{2}}(\V{a}_i^K - \V{a}_k^K)$, respectively.
Thus, the angle between edges $(\V{a}_j^K - \V{a}_k^K)$ and $(\V{a}_i^K - \V{a}_k^K)$ in the $\mathbb{D}_K^{-1}$-norm
and that between $\V{q}_i^K$ and $\V{q}_j^K$ in the $\mathbb{D}_K$ norm are related by
\beq
\arccos \left (\frac{(\V{q}_i^K)^T \mathbb{D}_K \V{q}_j^K}
{\| \V{q}_i^K\|_{\mathbb{D}_K} \; \| \V{q}_j^K\|_{\mathbb{D}_K}} \right )
+ \arccos \left (\frac{(\V{a}_i^K-\V{a}_k^K)^T \mathbb{D}_K^{-1} (\V{a}_j^K-\V{a}_k^K)}
{\| (\V{a}_i^K-\V{a}_k^K)\|_{\mathbb{D}_K^{-1}} \; \| (\V{a}_j^K-\V{a}_k^K)\|_{\mathbb{D}_K^{-1}}} \right )
= \pi .
\label{angles-1}
\eeq
Since (\ref{g-nonobtuse}) means the first angle on the left-hand side of the above equation is between $\pi/2$ and
$\pi$, we conclude that {\em condition (\ref{g-nonobtuse}) is equivalent to the requirement that
the angles of elements be non-obtuse when measured in the $\mathbb{D}_K^{-1}$ norm.}

It should be emphasized that condition (\ref{g-nonobtuse}) has been obtained by considering
only local stiffness matrices on individual elements. For the current 2D situation, this means that
each term in (\ref{stiffness-matrix-1}) has been required to be non-positive. Clearly, this is too strong
since we only need $a_{ij} \le 0$ for $i\neq j$ for $A$ to be an $M$-matrix.
For the special case $\mathbb{D} = I$, the Delaunay condition requiring
the sum of any pair of angles opposite a common edge
to be less than or equal to $\pi$ (cf. Fig. \ref{f1}) is sufficient to guarantee $a_{ij} \le 0$ for $i\neq j$.
It is then natural to ask if condition (\ref{g-nonobtuse}) can be weakened and a Delaunay-type condition
exists for the general diffusion matrix $\mathbb{D}$.
This issue is studied in the next section.

% section 4
\section{A Delaunay-type mesh condition}

In this section, we develop a Delaunay-type mesh condition under
which the linear finite element scheme (\ref{disc-3}) satisfies DMP.  The main result is given in Theorem \ref{thm4.1}.
Its proof is broken into a series of Lemmas.

\begin{lem}
\label{lem4.1}
For any element $K$,
\beq
|K| (\nabla \phi_i|_K )^T \nabla \phi_j|_K  = -\frac{1}{2} \cot (\alpha_{ij}^K ), \quad i \neq j, \; i,j = 1, 2, 3
\label{lem4.1-1}
\eeq
where $\alpha_{ij}^K $ is the angle between edges $e_{ki}$ and $e_{kj}$, with $\V{a}_k^K$ being the third vertex.
\end{lem}

{\bf Proof.} This result has been obtained in \cite{For89}. For completeness, we give a short proof here.
Without loss of generality, we consider the case with $i=1$, $j = 2$, and $k=3$ (cf. Fig. \ref{f2}).
From (\ref{dihedral}), (\ref{q-phi}), and (\ref{q-h}),  we have
\bey
|K| (\nabla \phi_1|_K )^T \nabla \phi_2|_K & = & |K| (\V{q}_1^K)^T \V{q}_2^K
\nn \\
& = & |K|\; \| \V{q}_1^K\| \cdot \| \V{q}_2^K\| \cos(\pi - \alpha_{12}^K) 
\nn \\
& = & - \frac{|K|}{h_1^K h_2^K} \cos( \alpha_{12}^K) .
\nn
\eey
From Fig. \ref{f2}, it is easy to see
\[
|K| = \frac{1}{2} h_2^K \|\V{a}_1^K - \V{a}_3^K \| = \frac{h_1^K h_2^K}{2 \sin(\alpha_{12}^K)} .
\]
Combining the above results, we obtain inequality (\ref{lem4.1-1}).
\proofend

\vspace{10pt}

The angle $\alpha_{ij}^K$ can be calculated  in terms of the $\V{q}$ vectors as in (\ref{dihedral})
or in terms of the edge vectors as
\beq
\alpha_{ij}^K = \arccos \left ( \frac{ (\V{a}_i^K - \V{a}_k^K)^T (\V{a}_j^K - \V{a}_k^K)}{\| \V{a}_i^K - \V{a}_k^K\| \;
\|\V{a}_j^K - \V{a}_k^K\| } \right ) .
\label{dihedral-1}
\eeq
The above formula is more desirable if linear coordinate transformations are involved.
This is because,
under a linear coordinate transformation, the edge vectors of $K$ will remain to be
the edge vectors of the transformed element but in general the $\V{q}$ vectors will not.
The latter is due to the fact
that orthogonality between vectors is not preserved by linear coordinate transformations. 

\begin{lem}
\label{lem4.2}
For any element $K$,
\beq
|K| (\nabla \phi_i|_K)^T \mathbb{D}_K \nabla \phi_j|_K =
-\frac{\sqrt{\mbox{\em det}(\mathbb{D}_K)}}{2} \cot (\alpha_{ij, \mathbb{D}_K^{-1}}^K ), \quad i \neq j, \; i,j = 1, 2, 3
\label{lem4.2-1}
\eeq
where $\alpha_{ij,\mathbb{D}_K^{-1}}^K$ is the angle between edges $e_{ki}$ and $e_{kj}$
(with $\V{a}_k^K$ being the third vertex)
measured in the metric $\mathbb{D}_K^{-1}$, i.e.,
\beq
\alpha_{ij,\mathbb{D}_K^{-1}}^K = \arccos \left ( \frac{ (\V{a}_i^K - \V{a}_k^K)^T \mathbb{D}_K^{-1} (\V{a}_j^K - \V{a}_k^K)}
{\| \V{a}_i^K - \V{a}_k^K\|_{\mathbb{D}_K^{-1}} \; \|\V{a}_j^K - \V{a}_k^K\|_{\mathbb{D}_K^{-1}} } \right ) .
\label{dihedral-2}
\eeq
\end{lem}

{\bf Proof.} Consider a linear mapping $G: K \to \tilde{K}$ defined as
\beq
\left (\begin{array}{c} \xi \\ \eta \end{array} \right ) = \mathbb{D}_K^{-\frac{1}{2}}
\left (\begin{array}{c} x \\ y \end{array} \right ) ,\quad \forall \; (x,y) \in K
\label{linear-map}
\eeq
where $\tilde{K} = G(K)$ and $(x,y)$ and $(\xi, \eta)$ are the coordinates in $K$ and $\tilde{K}$,
respectively. Let $\tilde{\V{a}}_i = G(\V{a}_i^K)$ ($i = 1, 2, 3$), $\tilde{e}_{ij} = G(e_{ij})$ ($i \neq j$),
and $\tilde{\nabla} = ((\partial/\partial \xi), (\partial/\partial \eta))^T$. Denote the angles of $\tilde{K}$
by $\tilde{\alpha}_{ij}$. It is easy to show that $\tilde{e}_{ij}$'s form the edges of $\tilde{K}$
and $\tilde{\phi}_i(\xi, \eta) \equiv \phi_i|_K (F^{-1}(\xi,\eta))$ ($i=1, 2, 3$) form the linear basis functions on $\tilde{K}$.
Moreover, 
\[
\nabla = \mathbb{D}_K^{-\frac{1}{2}} \tilde{\nabla } .
\] 

Since $\nabla \phi_i $ and $\nabla \phi_j$ are constant on $K$, we have
\bey
|K| (\nabla \phi_i|_K)^T \mathbb{D}_K \nabla \phi_j|_K
& = & \int_K (\nabla \phi_i)^T \mathbb{D}_K \nabla \phi_j d x d y
\nn \\
& = & \int_{\tilde{K}} (\tilde{\nabla} \tilde{\phi}_i)^T \tilde{\nabla} \tilde{\phi}_j
\mbox{det}(\mathbb{D}_K^{\frac{1}{2}} )d \xi d \eta
\nn \\
& = &  \sqrt{\mbox{det}(\mathbb{D}_K)} \; |\tilde{K}| \;
(\tilde{\nabla} \tilde{\phi}_i|_{\tilde{K}})^T \tilde{\nabla} \tilde{\phi}_j |_{\tilde{K}}.
\nn
\eey
Applying Lemma \ref{lem4.1} to the last term in the above equation on element $\tilde{K}$, we have
\beq
|K| (\nabla \phi_i|_K )^T \mathbb{D}_K \nabla \phi_j|_K
= -\frac{\sqrt{\mbox{det}(\mathbb{D}_K)}}{2} \cot (\tilde{\alpha}_{ij} ) .
\label{lem4.2-3}
\eeq
From 
\[
\tilde{\V{a}}_i - \tilde{\V{a}}_k = \mathbb{D}_K^{-\frac{1}{2}} (\V{a}_i^K - \V{a}_k^K),\quad
\| \tilde{\V{a}}_i - \tilde{\V{a}}_k \| = \| \V{a}_i^K - \V{a}_k^K \|_{\mathbb{D}_K^{-1}}
\]
and similar formulas for $(\V{a}_j^K - \V{a}_k^K)$, $\tilde{\alpha}_{ij}$ can be expressed as
\bey
\tilde{\alpha}_{ij} & = & \arccos \left ( \frac{ (\tilde{\V{a}}_i - \tilde{\V{a}}_k)^T (\tilde{\V{a}}_j - \tilde{\V{a}}_k)}
{\| \tilde{\V{a}}_i - \tilde{\V{a}}_k\| \; \|\tilde{\V{a}}_j - \tilde{\V{a}}_k\| } \right ) 
\nn \\
& = & \arccos \left ( \frac{ (\V{a}_i^K - \V{a}_k^K)^T \mathbb{D}_K^{-1} (\V{a}_j^K - \V{a}_k^K)}
{\| \V{a}_i^K - \V{a}_k^K\|_{\mathbb{D}_K^{-1}} \; \|\V{a}_j^K - \V{a}_k^K\|_{\mathbb{D}_K^{-1}} } \right )
\nn
\\
& = & \alpha_{ij,\mathbb{D}_K^{-1}}^K .
\nn
\eey
Combining this result with (\ref{lem4.2-3}) gives (\ref{lem4.2-1}).
\proofend

\begin{lem}
\label{lem4.3}
The entry $a_{ij}$ of
the stiffness matrix $A$, (\ref{stiffness-matrix-1}), can be expressed as
\beq
a_{ij} = -\frac{\sqrt{\mbox{\em det}(\mathbb{D}_K)}}{2} \cot (\alpha_{ij, \mathbb{D}_K^{-1}}^K )
-\frac{\sqrt{\mbox{\em det}(\mathbb{D}_{K'})}}{2} \cot (\alpha_{ij, \mathbb{D}_{K'}^{-1}}^{K'} ) .
\label{lem4.3-1}
\eeq
\end{lem}

{\bf Proof.} This lemma is a consequence of combination of (\ref{stiffness-matrix-1}) and Lemma \ref{lem4.2}.
\proofend

\begin{thm}
\label{thm4.1}
If the triangular mesh satisfies
\bey
&& \frac{1}{2}\left [
\alpha_{ij, \mathbb{D}_{K}^{-1}}^K + \alpha_{ij, \mathbb{D}_{K'}^{-1}}^{K'}
+  \mbox{\em arccot}\left (\sqrt{\frac{\mbox{\em det}(\mathbb{D}_{K})}{\mbox{\em det}(\mathbb{D}_{K'})}}
\cot (\alpha_{ij, \mathbb{D}_{K}^{-1}}^K ) \right )
\right .
\nn \\
&& \qquad \qquad \left .
+  \; \mbox{\em arccot}\left (\sqrt{\frac{\mbox{\em det}(\mathbb{D}_{K'})}{\mbox{\em det}(\mathbb{D}_{K})}}
\cot (\alpha_{ij, \mathbb{D}_{K'}^{-1}}^{K'} ) \right ) \right ]
\le \pi ,  \quad \mbox{ for all interior edges $e_{ij}$}
\label{thm4.1-1}
\eey
where $K$ and $K'$ are the elements sharing $e_{ij}$, then
the linear finite element scheme (\ref{disc-3}) satisfies DMP.
\end{thm}

{\bf Proof.} We first show that if the mesh satisfies 
\beq
\alpha_{ij, \mathbb{D}_K^{-1}}^K + \mbox{arccot}
\left (\sqrt{\frac{\mbox{det}(\mathbb{D}_{K'})}{\mbox{det}(\mathbb{D}_{K})}}
\cot (\alpha_{ij, \mathbb{D}_{K'}^{-1}}^{K'} ) \right )
\le \pi, \quad \mbox{ for all interior edges $e_{ij}$}
\label{thm4.1-1+1}
\eeq
then the conclusion holds. Indeed, notice that the inequality
\[
\frac{\sqrt{\mbox{det}(\mathbb{D}_K)}}{2} \cot (\alpha_{ij, \mathbb{D}_K^{-1}}^K )
+ \frac{\sqrt{\mbox{det}(\mathbb{D}_{K'})}}{2} \cot (\alpha_{ij, \mathbb{D}_{K'}^{-1}}^{K'} ) \ge 0
\]
can be written as
\[
\alpha_{ij, \mathbb{D}_K^{-1}}^K \le \mbox{arccot}\left ( - 
\sqrt{\frac{\mbox{det}(\mathbb{D}_{K'})}{\mbox{det}(\mathbb{D}_{K})}}
\cot (\alpha_{ij, \mathbb{D}_{K'}^{-1}}^{K'} ) \right )
=\pi - \mbox{arccot}\left ( 
\sqrt{\frac{\mbox{det}(\mathbb{D}_{K'})}{\mbox{det}(\mathbb{D}_{K})}}
\cot (\alpha_{ij, \mathbb{D}_{K'}^{-1}}^{K'} ) \right ),
\]
which is exactly (\ref{thm4.1-1+1}). Then, from Lemma \ref{lem4.3} we have
$a_{ij} \le 0$ for $i = 1, ..., N_{vi}$ and $j = 1, ..., N_v$ if (\ref{thm4.1-1+1}) is satisfied.
The result also means $a_{ij} \le 0$ for all $i\neq j$ due to the special structure  (\ref{fem-matrix}) of the
stiffness matrix.
Following the proof of Theorem 2.1 of \cite{LH10} we can then show that $A$ is an $M$-matrix and
has non-negative row sums, which implies that the linear finite element scheme (\ref{disc-3}) satisfies
DMP (cf. Stoyan \cite{Sto86} or Lemma 1.2 of \cite{LH10}).

Next, it is easy to show that (\ref{thm4.1-1+1}) is equivalent to
\beq
\alpha_{ij, \mathbb{D}_{K'}^{-1}}^{K'} + \mbox{arccot}
\left (\sqrt{\frac{\mbox{det}(\mathbb{D}_{K})}{\mbox{det}(\mathbb{D}_{K'})}}
\cot (\alpha_{ij, \mathbb{D}_{K}^{-1}}^K ) \right ) 
\le \pi  .
\label{thm4.1-1+2}
\eeq
As a result, (\ref{thm4.1-1}) and (\ref{thm4.1-1+1}) are mathematically equivalent.
\proofend

\vspace{12pt}

We now study the mesh condition (\ref{thm4.1-1}). We first consider the case with constant $\mathbb{D}$.
For this case,
\[
\mathbb{D}_{K} = \mathbb{D}_{K'} = \mathbb{D},\quad
\mbox{det}(\mathbb{D}_{K}) = \mbox{det}(\mathbb{D}_{K'}) = \mbox{det}(\mathbb{D}) .
\]
Then (\ref{thm4.1-1}) reduces to
\beq
\alpha_{ij, \mathbb{D}^{-1}}^{K} + \alpha_{ij, \mathbb{D}^{-1}}^{K'} \le \pi,
\quad \mbox{ for all interior edges $e_{ij}$} .
\label{thm4.1-2}
\eeq
For the special case with $\mathbb{D} = I$, (\ref{thm4.1-2}) reduces to
\beq
\alpha_{ij}^{K} + \alpha_{ij}^{K'} \le \pi,
\quad \mbox{ for all interior edges $e_{ij}$}
\label{delaunay-cond}
\eeq
which is exactly the Delaunay condition (cf. Fig. \ref{f1}). Thus, 
{\em  the mesh condition (\ref{thm4.1-1}) reduces to the Delaunay condition for the special case
$\mathbb{D} = I$ and is a generalization of the Delaunay condition for a general $\mathbb{D}$.}

Next, we consider the mesh condition
\beq
\alpha_{ij, \mathbb{D}_K^{-1}}^K \le \frac{\pi}{2}, \quad i \neq j,\; i, j = 1, 2, 3,
\quad \forall \; K \in \mathcal{T}_h
\label{thm4.1-3}
\eeq
for a general matrix-valued function $\mathbb{D} = \mathbb{D}(x,y)$.
From (\ref{angles-1}) and (\ref{dihedral-2}) it is not difficult to see that this mesh condition
is equivalent to the anisotropic non-obtuse angle condition (\ref{g-nonobtuse}). Moreover,
under (\ref{thm4.1-3}) we have
\[
\alpha_{ij, \mathbb{D}_{K}^{-1}}^{K} \le \frac{\pi}{2},\quad
\cot (\alpha_{ij, \mathbb{D}_{K}^{-1}}^{K} ) \ge 0,\quad
\mbox{arccot}
\left (\sqrt{\frac{\mbox{det}(\mathbb{D}_{K})}{\mbox{det}(\mathbb{D}_{K'})}}
\cot (\alpha_{ij, \mathbb{D}_{K}^{-1}}^{K} ) \right ) \le \frac{\pi}{2},
\]
and similar results for $\alpha_{ij, \mathbb{D}_{K'}^{-1}}^{K'}$ and thus 
(\ref{thm4.1-1}) is true. Therefore, (\ref{thm4.1-3}), or equivalently (\ref{g-nonobtuse}),
implies (\ref{thm4.1-1}). In other words,
{\em the mesh condition (\ref{thm4.1-1}) is weaker than the anisotropic
non-obtuse angle condition (\ref{g-nonobtuse}).}

\begin{figure}[thb]
\centering
\hbox{
\hspace{10mm}
\begin{minipage}[b]{3in}
\centerline{(a): $\sqrt{\mbox{det}(\mathbb{D}_{K'})/\mbox{det}(\mathbb{D}_{K})} = 100$}
\includegraphics[width=3in]{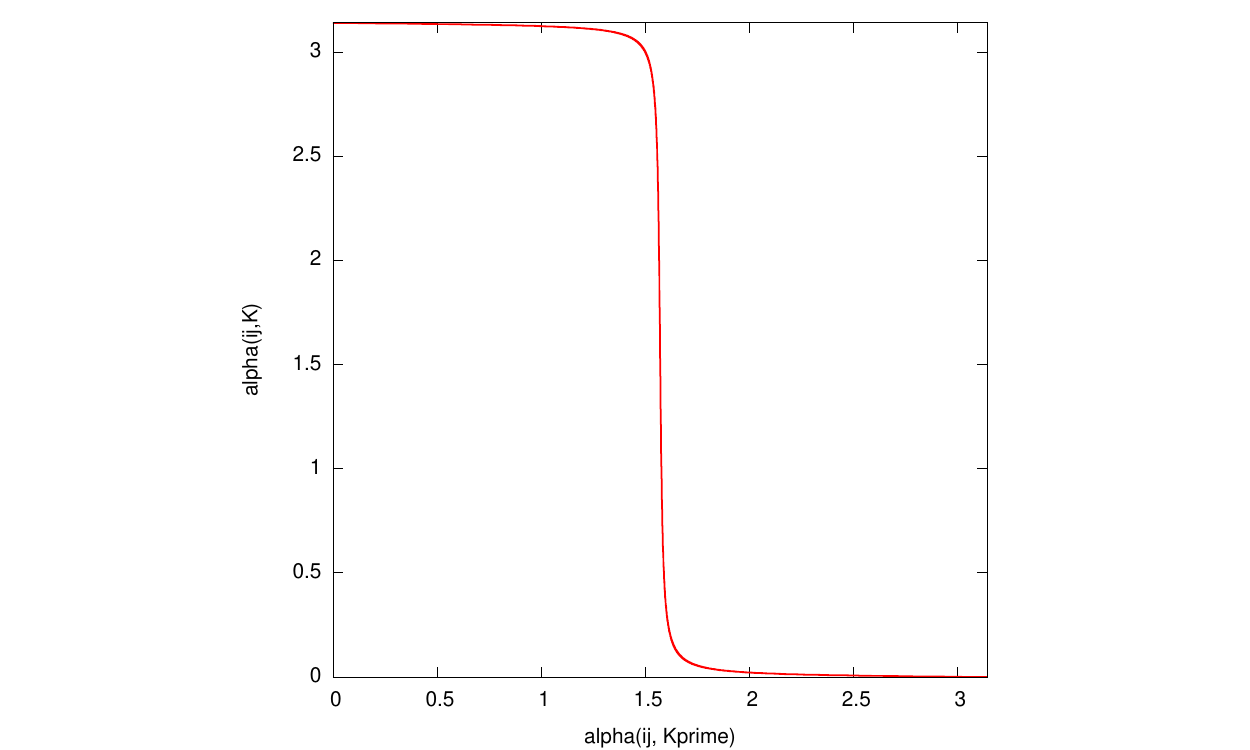}
\end{minipage}
\begin{minipage}[b]{3in}
\centerline{(b): $\sqrt{\mbox{det}(\mathbb{D}_{K'})/\mbox{det}(\mathbb{D}_{K})} = 0.01$}
\includegraphics[width=3in]{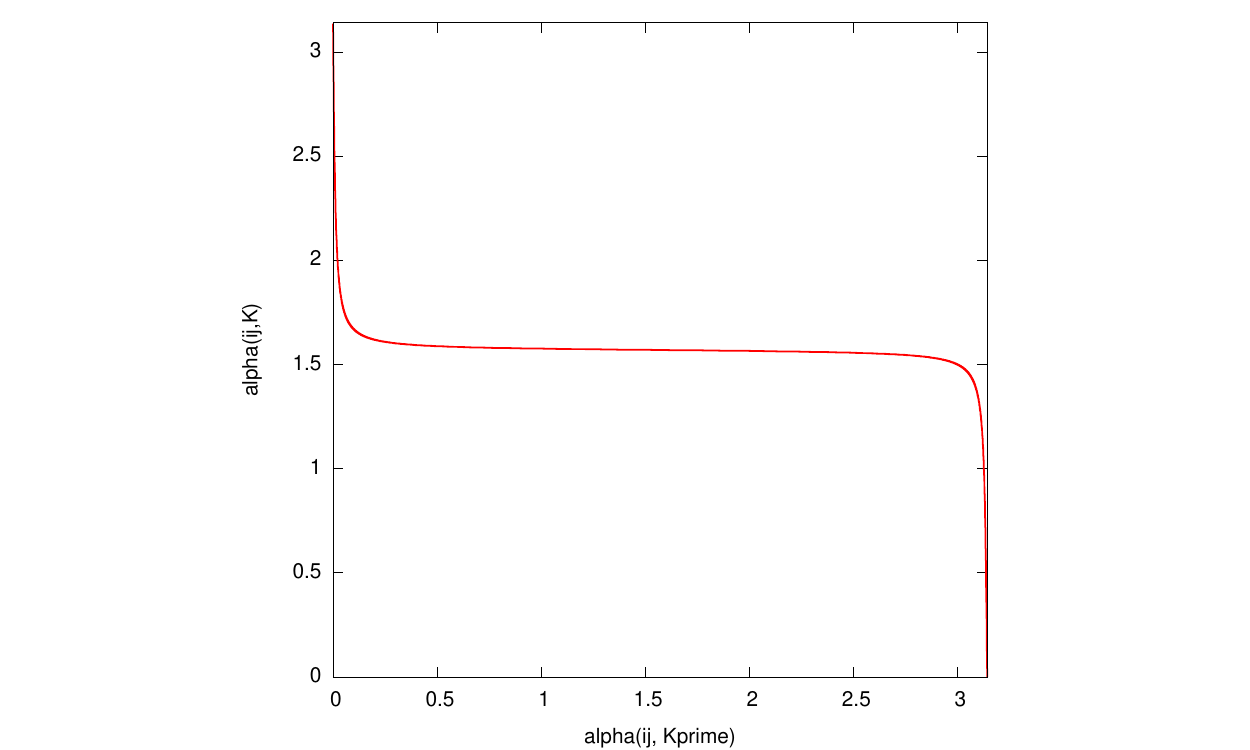}
\end{minipage}
}
\caption{The $x$ and $y$ axes are $\alpha_{ij, \mathbb{D}_{K'}^{-1}}^{K'}$ and $\alpha_{ij, \mathbb{D}_K^{-1}}^K$,
respectively. The DMP satisfaction region (satisfying the mesh condition (\ref{thm4.1-1})) is below
the plotted curve.
}
\label{f3}
\end{figure}

Finally, we consider some special cases for (\ref{thm4.1-1}).
It is obvious that (\ref{thm4.1-1}) reduces to (\ref{thm4.1-2})
when $\mbox{det}(\mathbb{D}_{K}) = \mbox{det}(\mathbb{D}_{K'})$.
In Fig. \ref{f3}  the region of ($\alpha_{ij, \mathbb{D}_{K'}^{-1}}^{K'}, \alpha_{ij, \mathbb{D}_K^{-1}}^K$)
satisfying the mesh condition (\ref{thm4.1-1}) is plotted for two cases where
the ratio $\mbox{det}(\mathbb{D}_{K'})/\mbox{det}(\mathbb{D}_{K})$ is either large or small.
From the figure, one can see that when the ratio is large (Fig. \ref{f3}(a)),
$\alpha_{ij, \mathbb{D}_{K'}^{-1}}^{K'}$ should essentially be non-obtuse
whereas $\alpha_{ij, \mathbb{D}_{K}^{-1}}^K$ can basically be any angle between $0$ and $\pi$.
On the other hand, when the ratio is small (Fig. \ref{f3}(b)), the roles of
$\alpha_{ij, \mathbb{D}_{K}^{-1}}^{K}$ and $\alpha_{ij, \mathbb{D}_{K'}^{-1}}^{K'}$ switch.
This observation is consistent with that made by Xu and Zikatanov \cite{XZ99}
that the non-obtuse angle condition should be imposed at edges where the diffusion coefficient
is discontinuous (and thus the ratio $\mbox{det}(\mathbb{D}_{K'})/\mbox{det}(\mathbb{D}_{K})$ can be
large or small) to guarantee DMP satisfaction. It is also interesting to observe from Fig. \ref{f3} that
the DMP satisfaction region overlaps with
$\alpha_{ij, \mathbb{D}_K^{-1}}^K + \alpha_{ij, \mathbb{D}_{K'}^{-1}}^{K'} \ge \pi$.

% f3

% section 5
\section{Numerical results}
\label{SEC:numerical}

In this section, we present some numerical results obtained for
BVP (\ref{bvp-pde}) and (\ref{bvp-bc}) with
\[
f \equiv 0,\quad g(x,0)=g(16,y)=0, 
\]
\beq
\nn 
g(0,y)= \left\{ \begin{array}{ll} 0.5y, & \text{ for } 0 \le y < 2 \\
 	1, & \text{ for  } 14 \le y \le 16 \end{array} \right. \text{ and  }\;
g(x,16)= \left\{ \begin{array}{ll} 1, & \text{ for  } 0 \le x \le 14 \\
 	8-0.5x, & \text{ for  } 14 < x \le  16. \end{array} \right.
\eeq
The diffusion matrix is defined as
\beq
\nn
\mathbb{D}(x,y) = \left ( \begin{array}{cc} 500.5 & 499.5 \\ 499.5 & 500.5 \end{array} \right ).
\eeq
This example has a constant but anisotropic diffusion matrix $\mathbb{D}$ and 
a continuous boundary condition. It satisfies the maximum principle and its solution
stays between $0$ and $1$.

\begin{figure}[thb]
\centering
\hbox{
\hspace{5mm}
\begin{minipage}[b]{3.5in}
\centerline{(a)}
\includegraphics[width=3.5in]{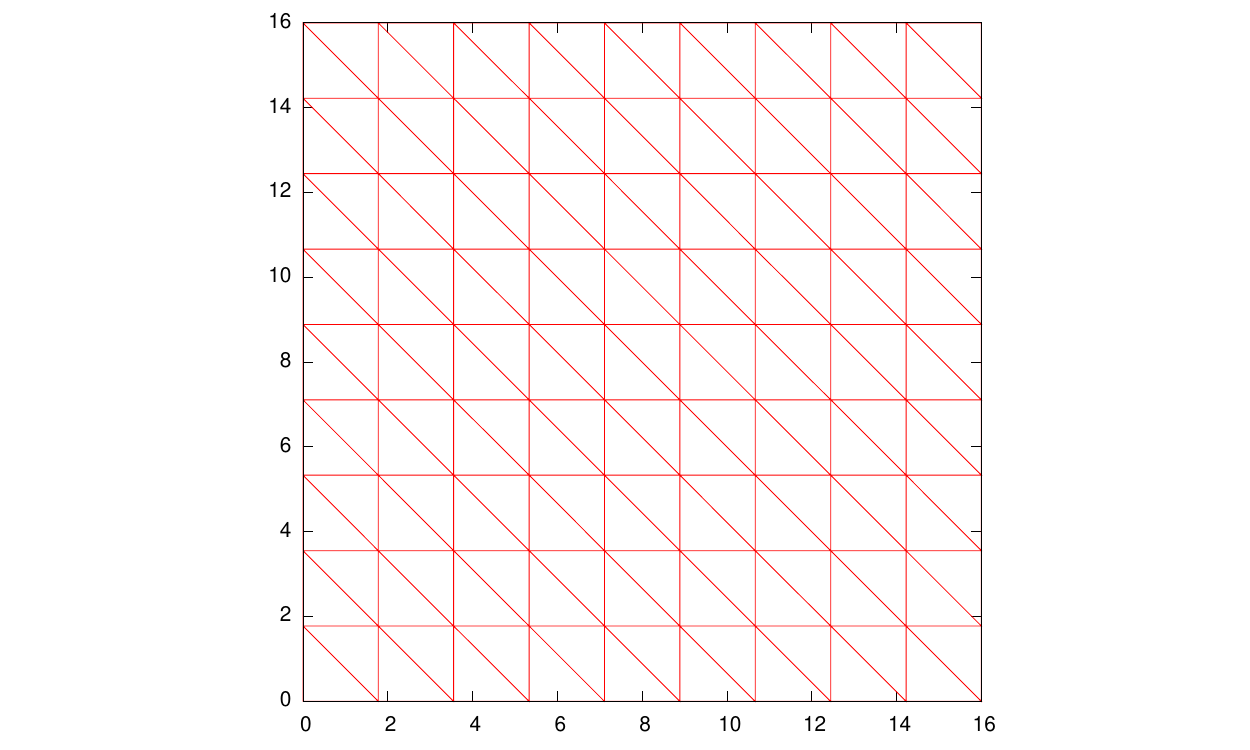}
\end{minipage}
\hspace{-25mm}
\begin{minipage}[b]{3.5in}
\centerline{(b)}
\includegraphics[width=3.5in]{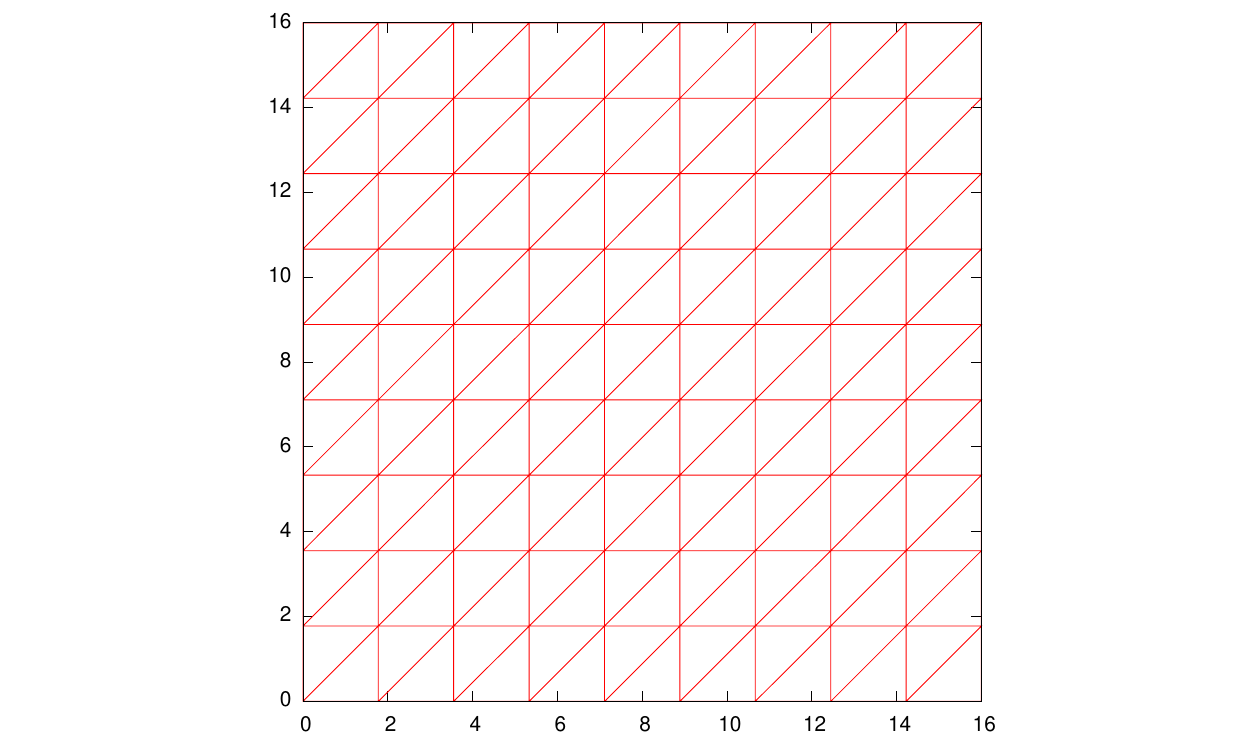}
\end{minipage}
}
\hbox{
\hspace{5mm}
\begin{minipage}[b]{3.5in}
\centerline{(c)}
\includegraphics[width=3.5in]{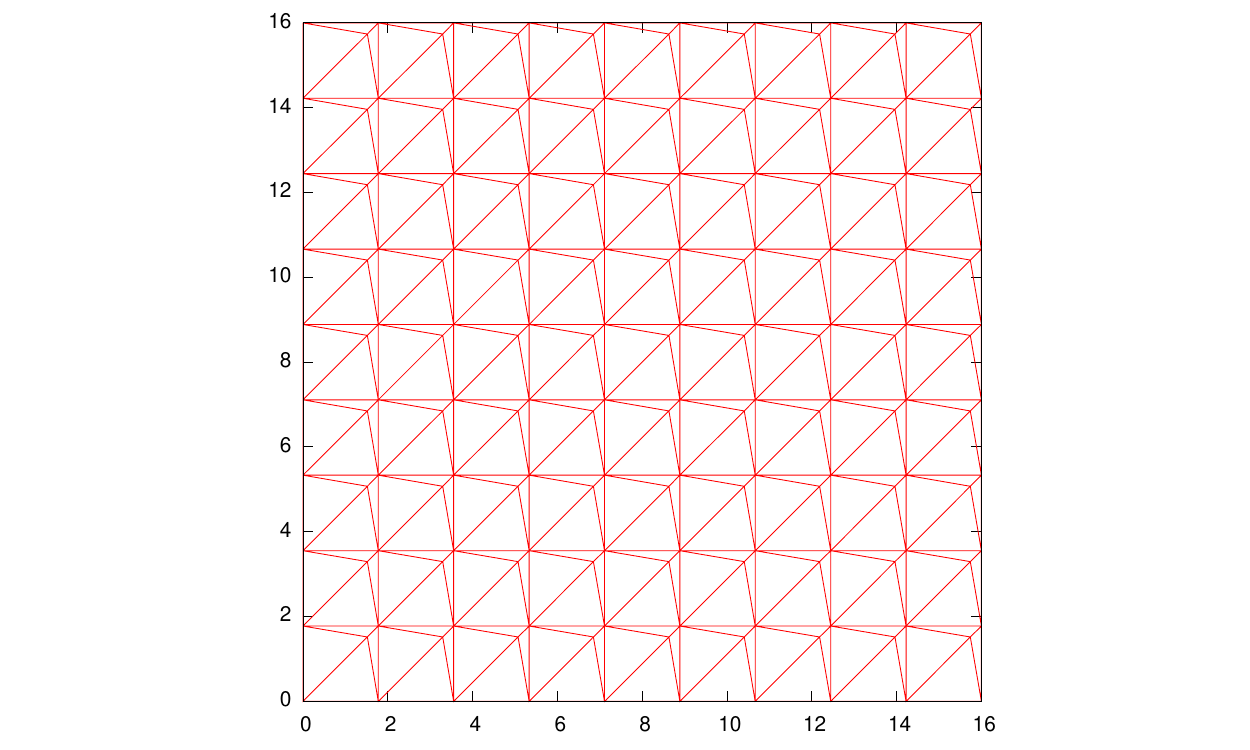}
\end{minipage}
\hspace{-25mm}
\begin{minipage}[b]{3.5in}
\centerline{(d)}
\includegraphics[width=3.5in]{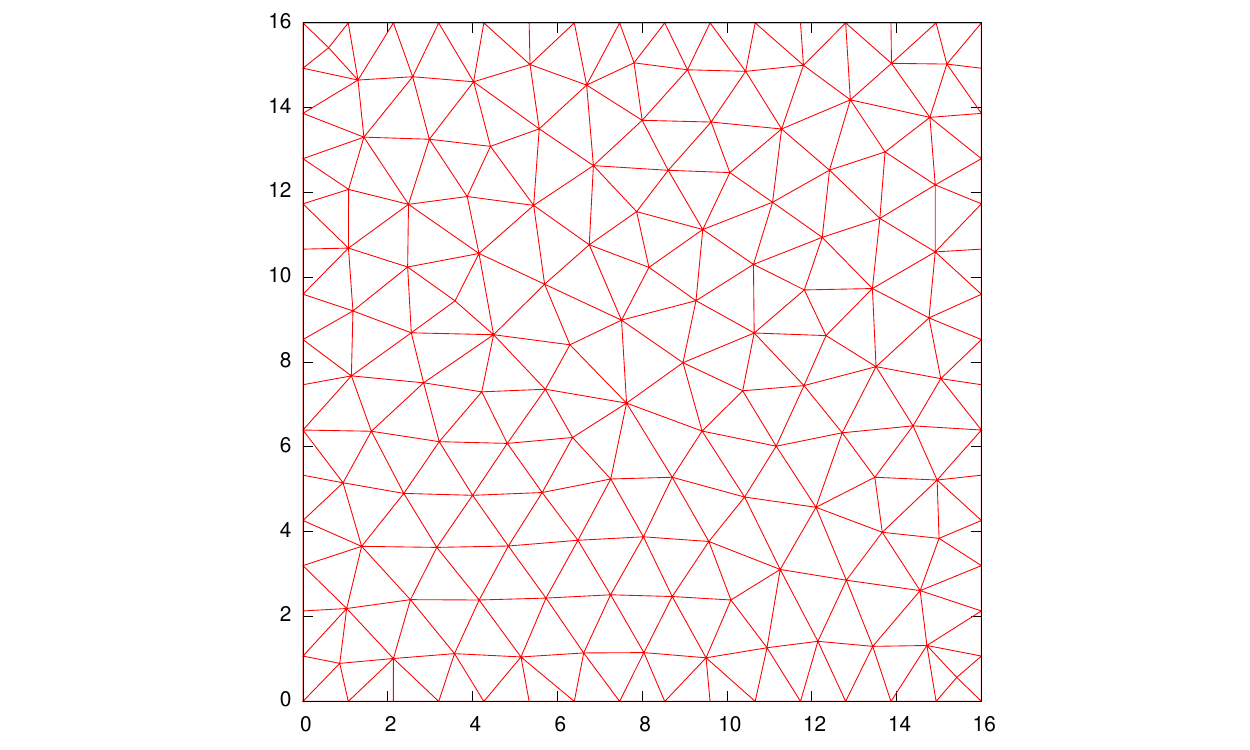}
\end{minipage}
}
\caption{Numerical example in Section \ref{SEC:numerical}. Meshes used in computation. The maximum values for
$\alpha_{ij, \mathbb{D}_K^{-1}}^K$ and $(\alpha_{ij, \mathbb{D}_K^{-1}}^K + \alpha_{ij, \mathbb{D}_{K'}^{-1}}^{K'})$,
respectively, are $0.98\pi$ and $1.96\pi$ for Mesh (a),  $0.49\pi$ and $0.98\pi$ for Mesh (b),
$0.51\pi$ and $\pi$ for Mesh (c), and  $0.98\pi$ and $1.96\pi$ for Mesh (d).}
\label{fex1-1}
\end{figure}

The computation is done with four types of triangular meshes shown in Fig. \ref{fex1-1}: Meshes (a) and
(b) are obtained by dividing a rectangle into two triangles using the northwest diagonal line and
the northeast line, respectively, Mesh (c) obtained by dividing a rectangle into four triangles with
the intersection toward the northeast corner, and Mesh (d) is a Delaunay mesh (which satisfies
the Delaunay condition). As mentioned in the previous section, mesh condition (\ref{thm4.1-1})
reduces to (\ref{thm4.1-2}) for the current example (with constant $\mathbb{D}$).
Note that Meshes (a) and (d) do not satisfy (\ref{thm4.1-2}) whereas Meshes (b) and (c) do (cf. Fig. \ref{fex1-1}).
Especially, Mesh (c) has obtuse elements (with angles greater than $\pi/2$ in the $\mathbb{D}^{-1}$--norm).

\begin{figure}[thb]
\centering
\hbox{
\hspace{25mm}
\begin{minipage}[b]{2in}
\centerline{(a)}
\includegraphics[width=2in]{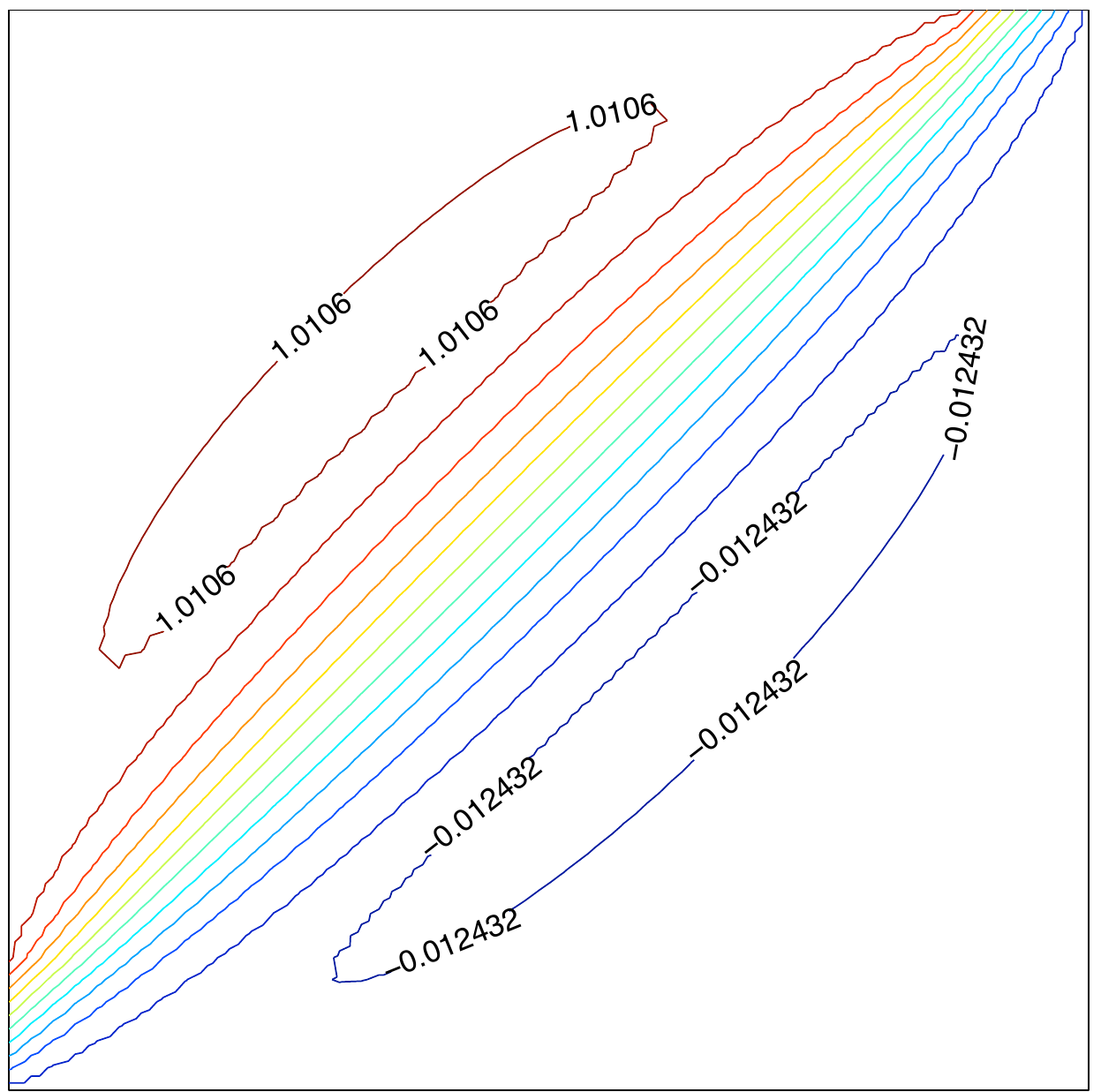}
\end{minipage}
\hspace{10mm}
\begin{minipage}[b]{2in}
\centerline{(b)}
\includegraphics[width=2in]{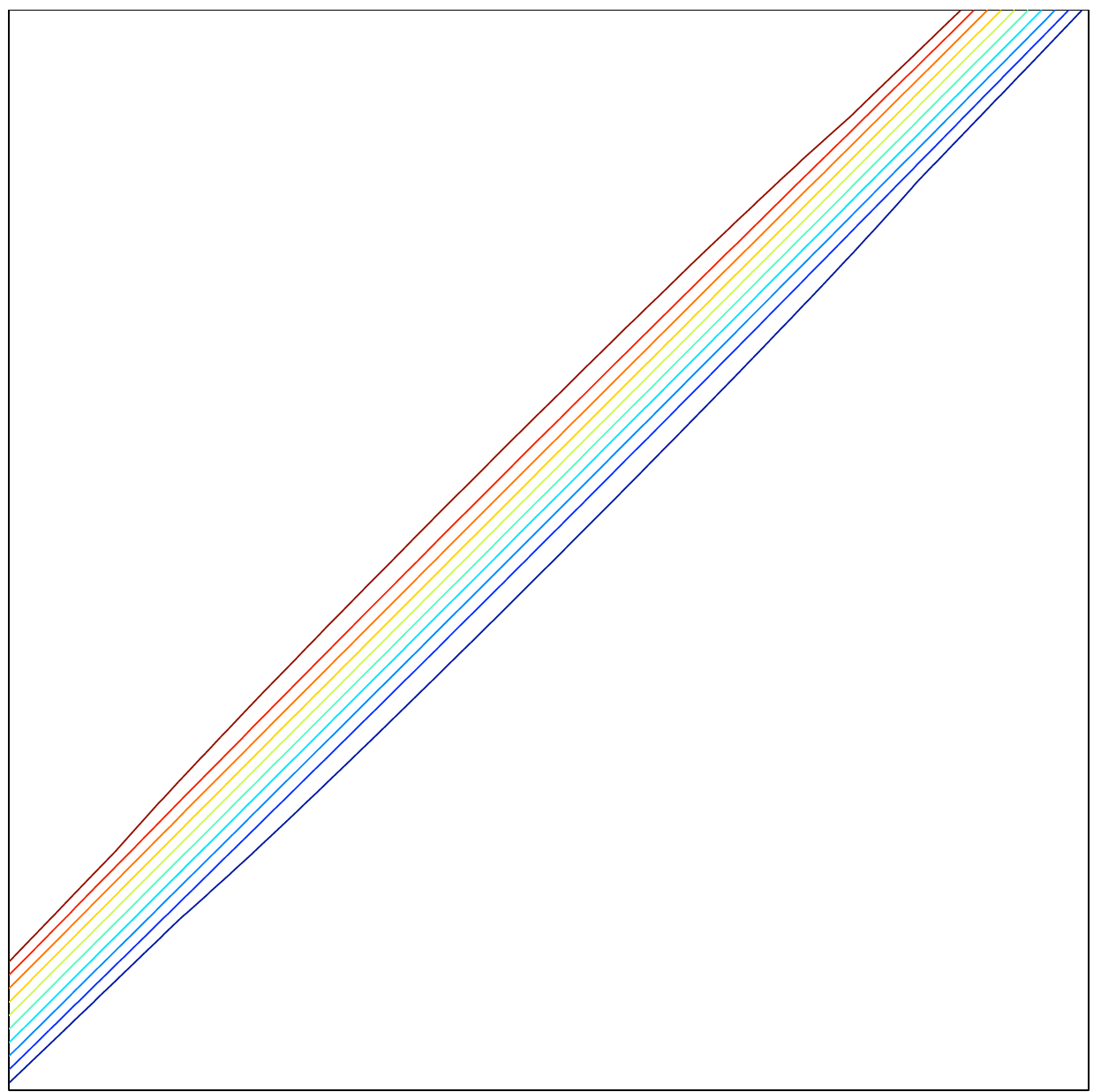}
\end{minipage}
}
\hbox{
\hspace{25mm}
\begin{minipage}[b]{2in}
\centerline{(c)}
\includegraphics[width=2in]{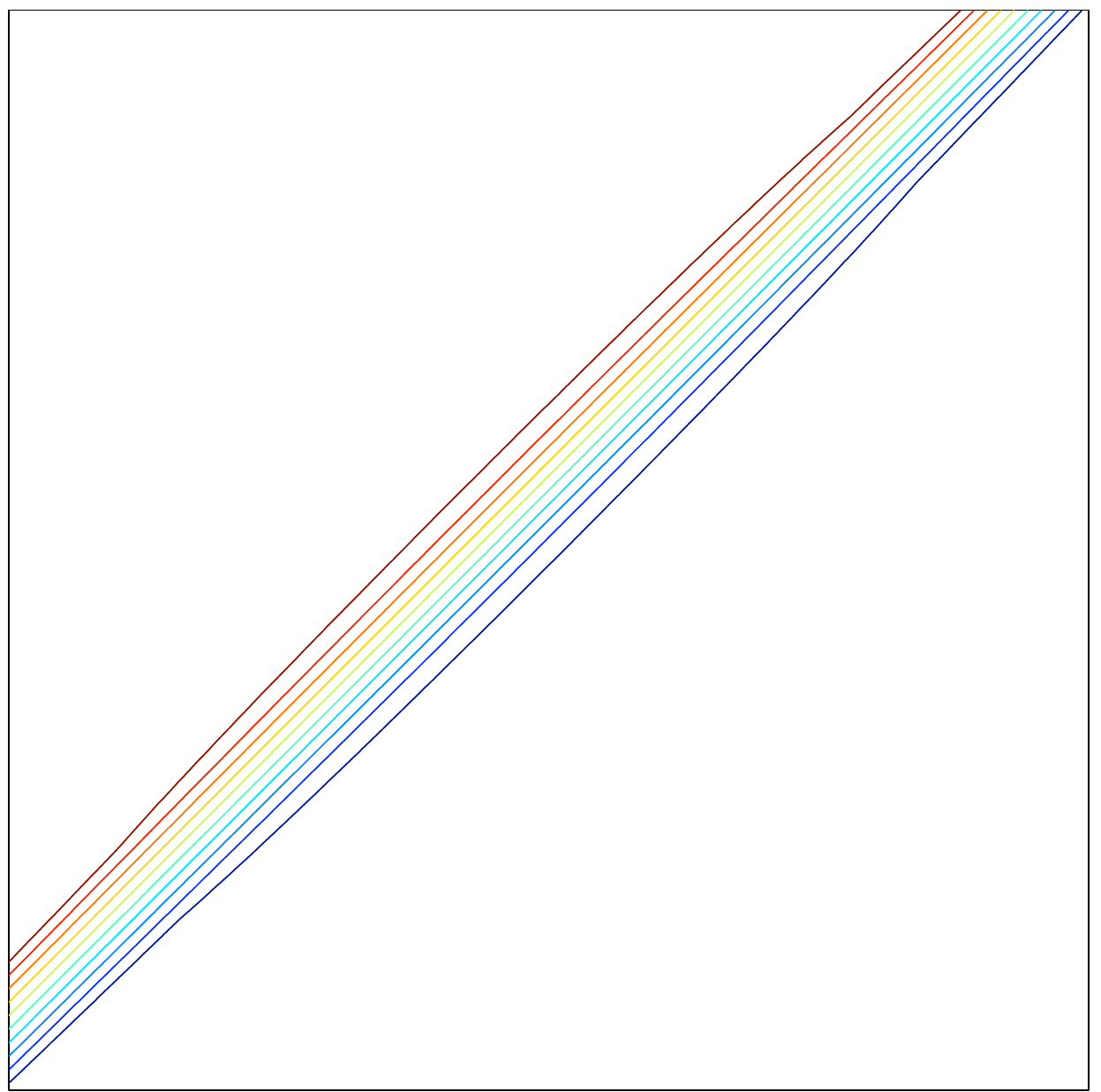}
\end{minipage}
\hspace{10mm}
\begin{minipage}[b]{2in}
\centerline{(d)}
\includegraphics[width=2in]{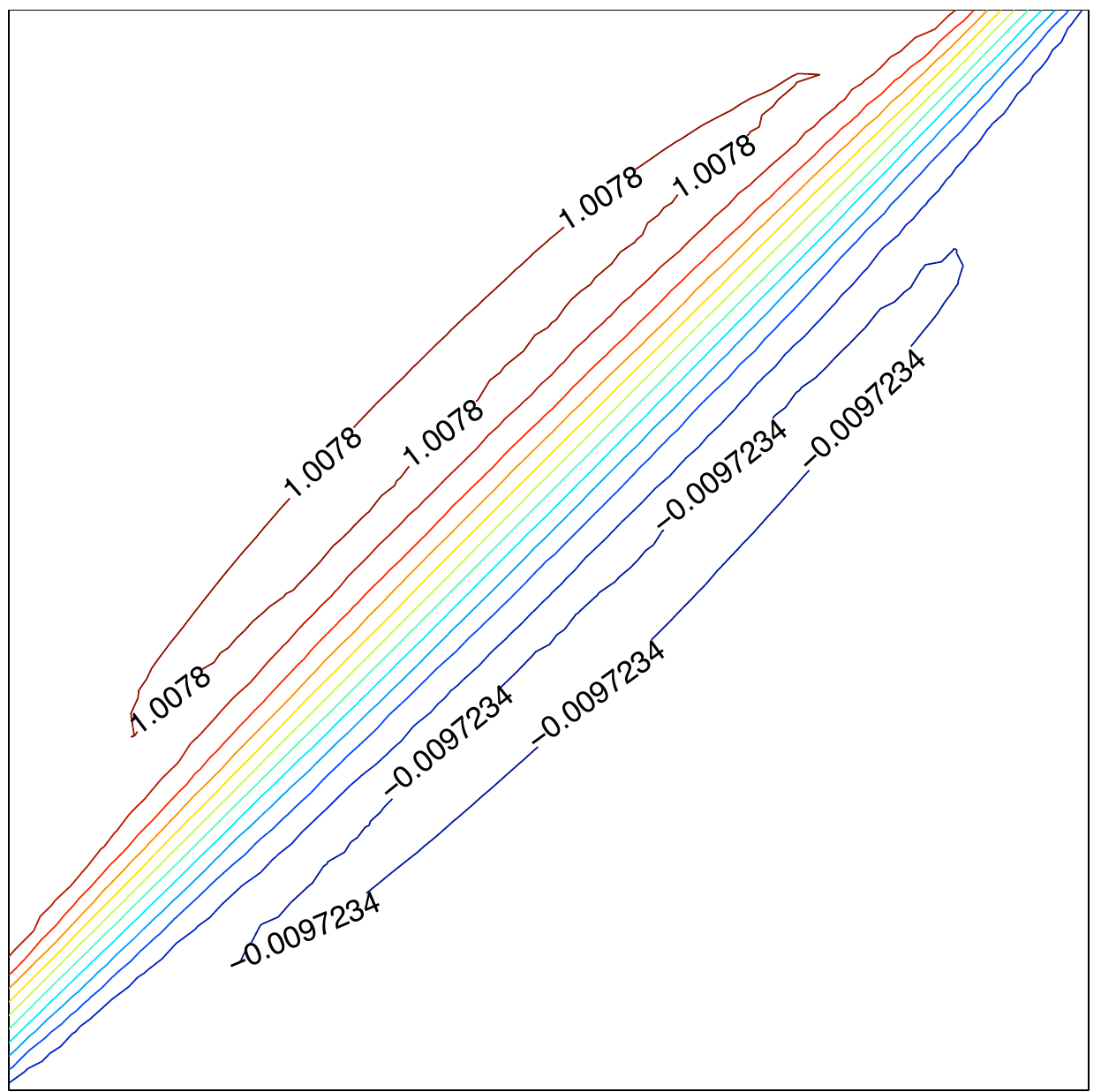}
\end{minipage}
}
\caption{Numerical example in Section \ref{SEC:numerical}. Contours of linear finite element solutions.}
% 4920 elements
\label{fex1-2}
\end{figure}

\begin{figure}[thb]
\centering
\hbox{
\hspace{2.5mm}
\begin{minipage}[b]{3in}
\centerline{(a): For the type of mesh in Fig. \ref{fex1-1}(a).}
\includegraphics[width=3in]{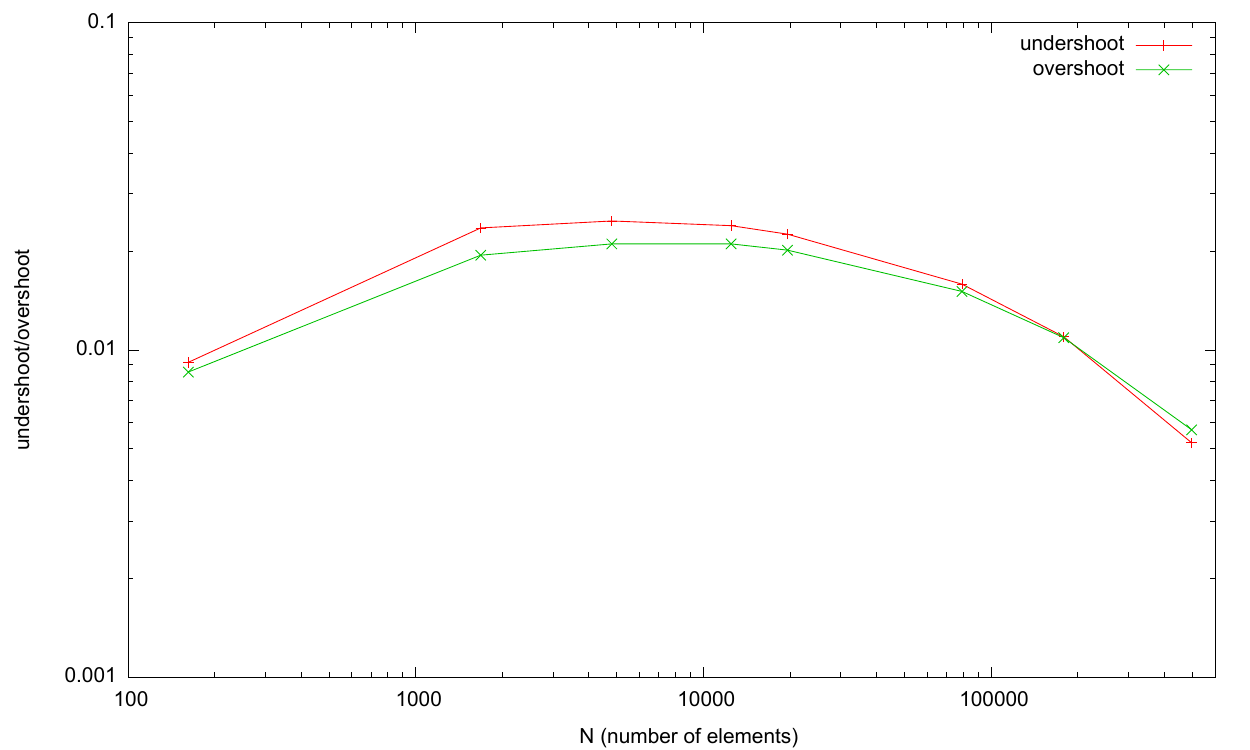}
\end{minipage}
\begin{minipage}[b]{3in}
\centerline{(b): For the type of mesh in Fig. \ref{fex1-1}(d).}
\includegraphics[width=3in]{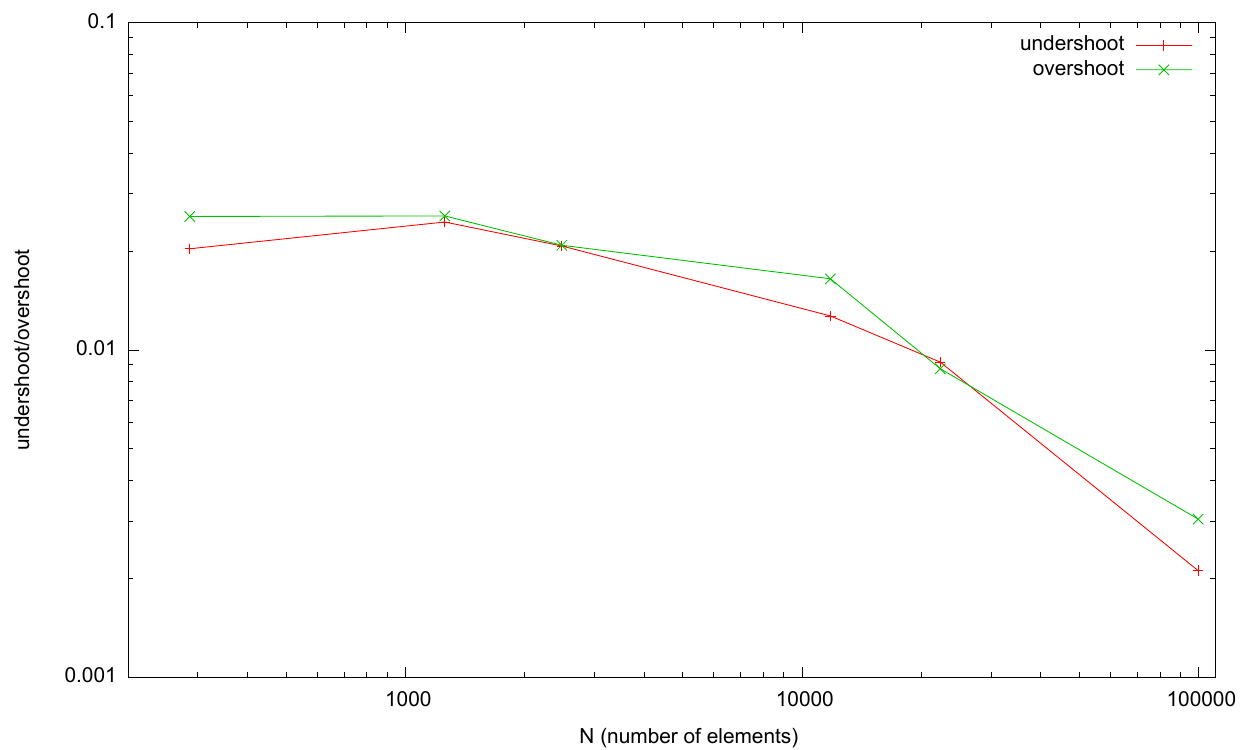}
\end{minipage}
}
\caption{Numerical example in Section \ref{SEC:numerical}. Overshoots and undershoots as functions of the number of
mesh elements.}
\label{fex1-3}
\end{figure}

Fig. \ref{fex1-2} shows the contours of the linear finite element solutions obtained for meshes finer than those
shown in Fig. \ref{fex1-1}. One can see that finite element solutions for both Meshes (b) and (c) stay between 0 and 1
and show no undershoots and overshoots. This is consistent with Theorem \ref{thm4.1}.
On the other hand, both Meshes (a) and (d) lead to undershoots and overshoots in the computational solutions.
Fig. \ref{fex1-3} shows these undershoots and overshoots as functions of the number of mesh elements.
As the mesh is refined, the undershoots and overshoots decrease very slowly and eventually reach
a rate $O(N^{-0.5})$, where $N$ is the number of elements.

% fex1-1

% fex1-2

% mesh (a): nbt, undershoot, overshoot
%162 9.17e-3 8.56e-3
%1682 2.36e-2 1.95e-2
%4802 2.48e-2 2.11e-2
%12482 2.40e-2 2.11e-2
%19602 2.26e-2 2.02e-2
%79202 1.59e-2 1.51e-2
%178802 1.096830e-02 1.0900e-2
%498002 5.201620e-03 5.700e-3

% mesh (d) (delaunay mesh)
% nbt, undershoot, overshoot
% 286 2.04e-2 2.56e-2
%1258 2.46e-2 2.57e-2
%2480 2.08e-2 2.09e-2
%11798 1.27e-2 1.65e-2
%22318 9.18e-3 8.76e-3
%99716 2.116050e-03 3.040e-3
%

% fex1-3

% section 6
\section{Conclusions and comments}

In the previous sections we have developed a Delaunay-type mesh condition (\ref{thm4.1-1}) under which
the linear finite element scheme (\ref{disc-3}) for solving the anisotropic diffusion problem
(\ref{bvp-pde}) and (\ref{bvp-bc}) satisfies DMP. This condition is weaker than the anisotropic
non-obtuse angle condition (\ref{g-nonobtuse}) developed in \cite{LH10}. It reduces to (\ref{thm4.1-2}) when
the diffusion matrix $\mathbb{D}$ is constant and especially to the Delaunay condition when $\mathbb{D} = I$.
The main theoretical result is given in Theorem \ref{thm4.1} and verified by numerical results.

It is well known that the Delaunay condition can be satisfied by a Delaunay mesh which
can be generated through edge swapping from an existing triangular mesh. Moreover,
Mlacnik and Durlofsky \cite{MD06} have demonstrated that a properly designed edge swapping procedure can improve
the monotonicity of finite volume approximations for anisotropic diffusion problems.
Clearly, the mesh condition (\ref{thm4.1-1}) can serve as a criterion for designing such a
procedure. The development of an edge swapping procedure based on (\ref{thm4.1-1}),
the convergence study of edge swapping, and the generation
of a mesh satisfying (\ref{thm4.1-1}) through edge swapping may deserve future investigation.
 
%\bibliographystyle{plain}
%\bibliography{/Users/huang/tex/bib/mmesh}

\end{document}